\documentclass[submission,copyright,creativecommons]{eptcs} 
\usepackage{amssymb}
\usepackage{breakurl}
\title{Towards a Paraconsistent Quantum Set Theory}
\author{Benjamin Eva}
\def\titlerunning{Towards a Paraconsistent Quantum Set Theory}
\def\authorrunning{B. E. Eva}
\begin{document}
\maketitle
\begin{abstract}
In this paper, we will attempt to establish a connection between quantum set theory, as developed by Ozawa, Takeuti and Titani (see, for example, \cite{Tit}, \cite{Tak}, \cite{Oza}), and topos quantum theory, as developed by Isham, Butterfield and D{\"o}ring, amongst others (see, for example, \cite{IB}, \cite{DI}, \cite{D1}). Towards this end, we will study algebraic valued set-theoretic structures whose truth values correspond to the clopen subobjects of the spectral presheaf of an orthomodular lattice of projections onto a given Hilbert space. In particular, we will attempt to recreate, in these new structures, Takeuti's original isomorphism between the set of all Dedekind real numbers in a suitably constructed model of set theory and the set of all self adjoint operators on a chosen Hilbert space.
\end{abstract}
\section{Introduction}

	We begin with some quick technical preliminaries. Given a complete Boolean algebra $B$, the `Boolean valued structure' $V^{(B)}$ is given by the following recursive definition, \\

	\hspace{2cm} $V_{\alpha}^{(B)} = \{x: func(x) \wedge ran(x) \subseteq B \wedge \exists  \xi <\alpha(dom(x) \subseteq V_{\xi}^{(B)})\}$\\ 

	\hspace{2cm}  $V^{(B)} = \{x: \exists \alpha (x \in V_{\alpha}^{(B)}) \}$\\ 

\hspace{2cm} i.e, $x \in V^{(B)}$ iff $fun(x) \wedge ran(x) \subseteq B \wedge dom(x) \subseteq V^{(B)}$ \\

	It is possible\footnote{For a full explanation of this and most other aspects of Boolean valued models, see \cite{Bell}} to treat $V^{(B)}$ as a new kind of model theoretic structure for ZF set theory, where any set theoretic sentence $\phi$ is given a truth value $\|\phi\|$ in the algebra $B$. Generally, $\|\phi\| \in \{0, 1\}$ won't hold. In the special cases where it does, we say that $\phi$ `holds/is true' or `fails/is false' in $V^{(B)}$ (depending on whether $\|\phi\| = 1$ or $\|\phi\| = 0$, respectively). It turns out that regardless of our choice of $B$, all of the axioms of ZFC are true in $V^{(B)}$. However, by varying $B$, we can obtain models that assign different truth values to sentences that can't be proved or disproved from the axioms of ZFC alone (this is one way to prove the independence of the continuum hypothesis). Intuitively, the idea is that each complete Boolean algebra $B$ corresponds to its own classical mathematical universe, many of which differ in very significant ways (truth of the continuum hypothesis, etc). 

	Now, since $V^{(B)}$ is a model of ZFC, there will always be an object $\mathbb{R} \in V^{(B)}$ such that $\|\mathbb{R}$ is the set of all Dedekind real numbers$\| = 1.$ We can then consider the set\footnote{It should be noted that this is a set in the sense of the `ground model' $V$, not in the sense $V^{(B)}$ itself.} $\mathbb{R}^{(B)} = \{u \in V^{(B)}| \|u \in \mathbb{R}\| = 1\}$ (this is called a `core' for $\mathbb{R}$). We can think of $\mathbb{R}^{(B)}$ as the set of all things that $V^{(B)}$ calls `Dedekind real numbers'. In \cite{Tak}, it was shown that if $B$ is an algebra of mutually compatible projection operators on some Hilbert space $H$, then $\mathbb{R}^{(B)}$ is isomorphic to the set of all self adjoint operators on $H$ whose spectral projections all lie in $B$. This was the founding result of quantum set theory (QST). Subsequently, there have been several attempts to extend this result to obtain, for any Hilbert space $H$, a corresponding algebraic valued model $V^{(A)}$ of some (not necessarily classical) set theory such that $\mathbb{R}^{(A)}$ is in bijective correspondence with the set $SA(H)$ of \emph{all} self adjoint operators on $H$. It has been argued (see, for example \cite{Oza 2}) that such a model could be useful for extending the usual Hilbert space formalism of quantum theory. Specifically, given such a model $V^{(A)}$, all physical propositions concerning the physical quantities of a given quantum system could be represented as propositions about the Dedekind reals in the model $V^{(A)}$. 

	This paper offers a new approach to this problem. In particular, we consider the possibility of using the complete Heyting algebra $Sub_{cl}(\underline{\Omega})$ (to be introduced in the next section) of clopen subobjects of the spectral presheaf on the projection lattice $P(H)$ of the relevant Hilbert space $H$ to build the structure $V^{Sub_{cl}(\underline{\Omega})}$. It is well known that for any complete Heyting algebra $H$, the structure $V^{({H})}$ will be a model of intuitionistic set theory. So this approach is situated in the context of intuitionistic rather than classical set theory. Some of the motivations of this project are as follows, \\

	\begin{itemize}
	\item{Until now, all the major attempts at finding a suitable model for the purposes of QST have used the orthomodular lattice $P(H)$ of projections on $H$ (see \cite{Oza}, \cite{Tit}). The non-distributivity of $P(H)$ has proved to be a significant and persistent obstacle in this regard. Specifically, $V^{(P(H))}$ is not a full model of any well known existing set theory. In topos quantum theory (TQT), $Sub_{cl}(\underline{\Omega})$ plays the role of the lattice of physical propositions associated with the quantum system represented by $H$ in the same way that $P(H)$ does in the orthodox Hilbert space formalism. It is hoped that the distributivity of $Sub_{cl}(\underline{\Omega})$ will significantly simplify some of the technical difficulties that beset the existing approaches to QST.}

	\item{On a conceptual level, QST shares a close affinity with TQT, in the sense that both approaches attempt to extend the Hilbert space formalism by employing non-standard mathematical universes with non-classical internal logics. It is hoped that studying the structure $V^{Sub_{cl}(\underline{\Omega})}$ will render the formal relationship between the two projects more perspicuous.}

	\item{In TQT, physical quantities are standardly represented by maps from the spectral presheaf to the `quantity valued object' in the relevant topos. One reason to think that this approach might be unsatisfactory is that the quantity valued object is not actually the same as the `real number object' in the topos. This seems incongruous in so far as TQT aims to utilise the internal structure of the relevant topos. The approach developed in this paper offers a new way to represent physical quantities in TQT (as Dedekind reals in $V^{Sub_{cl}(\underline{\Omega})}$) that has a closer connection to the ambient mathematical universe.}

	\end{itemize}

	\section{The Spectral Presheaf of an Orthomodular Lattice}

	In Topos Quantum Theory (TQT), we normally work with the spectral presheaf of the set of all bounded operators over some Hilbert space. Intuitively, this can be thought of as a way of `breaking a quantum system up into its constituent classical parts'. It is natural to ask whether a similar trick can be done for orthomodular lattices (which are often used to represent the logical structure of quantum theory). Specifically, it is natural to ask how much information we can get about an orthomodular lattice by looking just at its Boolean sub-algebras. It turns out that this question has already been given a fairly comprehensive answer by D{\"o}ring and Cannon \cite{Spec}. We will now give a quick overview of their main techniques and results.

		Firstly, recall that for any Boolean algebra $B$, the Stone space, $\Omega_{B}$ of $B$ is the set of all Boolean algebra homomorphisms from $B$ into $\{0, 1\}$. $\Omega_{B}$ can always be given a topology with basis sets of the form $U_{b} = \{\lambda \in \Omega_{B} | \lambda (b) = 1\}$ (where $b \in B$). This topology turns $\Omega_{B}$ into a compact totally disconnected Hausdorff space. Intuitively, Stone's theorem tells us that any Boolean algebra $B$ is isomorphic to the set $cl(\Omega_{B})$ of all clopen subsets of $\Omega_{B}$. Specifically, the map $\phi: B \rightarrow cl(\Omega_{B})$ defined by $\phi (b) = \{ \lambda \in \Omega_{B}| \lambda (b) = 1\}$ is always an isomorphism. Now, given an orthomodular lattice $L$, we can define the poset $B(L)$, of Boolean subalgebras of $L$, ordered by inclusion. This allows us to define the spectral presheaf of $L$, \\

	\textbf{Def 2.1:} The Spectral Presheaf, $\underline{\Omega}$ of an orthomodular lattice $L$ is the presheaf over $B(L)$ defined by \\

	\textbf{Objects:} Given $B \in B(L)$, $\underline{\Omega}_{B} = \Omega_{B}$\\

	\textbf{Arrows:} Given $B' \subseteq B$, $\underline{\Omega}_{B, B'}: \Omega_{B} \rightarrow \Omega_{B'}$\\

	\hspace{4.9 cm} $\lambda \mapsto \lambda|_{B'}$\\

	Now, it was shown by D{\"o}ring and Cannon that two orthomodular lattices are isomorphic (as orthomodular lattices) if and only if their spectral presheaves, as defined above, are isomorphic as presheaves. This tells us that the spectral presheaf $\underline{\Omega}$ is a complete invariant of the orthomodular lattice $L$. We now turn to the task of representing $L$ by clopen subobjects of $\underline{\Omega}$, thereby generalizing Stone's theorem to the new setting.

	First of all, given $a \in L$ and $B \in B(L)$, we define the daseinisation\footnote{Of course, we are assuming the completeness of $L$ here.} of $a$ at $B$ to be $\delta (a) = \bigwedge \{b \in B| b \geq a\}$.\\
	
	\textbf{Def 2.2:} Given $a \in L$, define the outer daseinisation presheaf $\underline{\delta(a)}$ over $B(L)$ by\\

	\textbf{Objects:} Given $B \in B(L)$, $\underline{\delta(a)}_{B} = \{\lambda \in \Omega_{B}| \lambda (\delta(a)) = 1\}$\\

	\textbf{Arrows:} Given $B' \subseteq B$, $\underline{\delta(a)}_{B, B'}: \delta(a)_{B} \rightarrow \delta(a)_{B'}$\\

	\hspace{5.8cm} $\lambda \mapsto \lambda|_{B'}$\\

	By Stone duality, $\underline{\delta(a)}_{B} \in cl(\Omega_{B})$, i.e. $\underline{\delta(a)}_{B}$ is a clopen subset of $\Omega_{B}$. Since this holds for all $B \in B(L)$, we call  $\underline{\delta(a)}$ a `clopen subobject' of $\underline{\Omega}$. It is easily shown that the lattice $Sub_{cl}(\underline{\Omega})$ of clopen subobjects of the spectral presheaf is a complete Heyting algebra. So, 2.2 defines a map $\underline{\delta}$ from $L$ into $Sub_{cl}(\underline{\Omega})$. It has been shown that this map satisfies the following properties, \\

	\textbf{Theorem 2.3:}\\

 (i) $\underline{\delta}$ is injective\\

(ii) $\underline{\delta}$ preserves all joins, i.e.  $\bigvee_{i \in I} \underline{\delta(a_{i})} = \underline{\delta( \bigvee_{i \in I} a_{i})}$, for any family \\ $\{a_{i}| i \in I\} \subseteq L$\\

(iii) $\underline{\delta(a)}$ is order preserving (monotone), i.e. $a \leq b$ in $L$ implies $\underline{\delta(a)} \leq \underline{\delta(b)}$ in $Sub_{cl}(\underline{\Omega})$\\

(iv) $\underline{\delta(0)} = \bot$, $\underline{\delta(1)} = \top$, where $0, 1$ are the minimal and maximal elements of $L$, respectively, and $\bot$ and $\top$ are the minimal and maximal elements of $Sub_{cl}(\underline{\Omega})$ ($\bot$ is the presheaf that takes each $B \in B(L)$ to the empty set and $\top$ is just $\underline{\Omega}$). \\

(v) $\underline{\delta(a \wedge b)} \leq \underline{\delta(a)} \wedge \underline{\delta(b)}$\\

	So we can think of $\underline{\delta}$ as an order preserving injection of $L$ into the complete Heyting algebra of clopen subobjects of the spectral presheaf of $L$. Since this map is monotone and join preserving, it has a monotone meet preserving upper adjoint, defined as follows, \\

\textbf{Def 2.4:} Given $\underline{S} \in Sub_{cl}(\underline{\Omega})$, define $\varepsilon (\underline{S}) = \bigvee \{a \in L| \underline{\delta(a)} \leq \underline{S}\}$ \\

	D{\"o}ring and Cannon showed that $\varepsilon$ has the following properties, \\

	\textbf{Theorem 2.5:}\\

	(i) $\varepsilon$ preserves all meets\\

	(ii) $\varepsilon$ is order preserving (monotone)\\

	(iii) $\varepsilon(\underline{\delta(a)}) = a$, for any $a \in L$.\\

	(iv) $\underline{\delta(\varepsilon(S))} \leq \underline{S}$, for any $\underline{S} \in Sub_{cl}(\underline{\Omega})$.\\

	(v) $\varepsilon(\underline{S} \vee \underline{T}) \geq \varepsilon(\underline{S}) \vee \varepsilon(\underline{T})$\\

	$\varepsilon$ can be used to define an equivalence relation on $Sub_{cl}(\underline{\Omega})$, defined by $\underline{S} \sim \underline{T}$ iff $\varepsilon (\underline{S}) = \varepsilon(\underline{T})$. We let $E$ denote the set of all equivalence classes of $Sub_{cl}(\underline{\Omega})$ under this equivalence relation. $E$ can be turned into a complete lattice by defining $\bigwedge_{i \in I} [\underline{S}_{i}] = [\bigwedge_{i \in I} \underline{S}_{i}]$, $[\underline{S}] \leq [\underline{T}]$ iff $[\underline{S}] \wedge [\underline{T}] =  [\underline{S}]$ and $\bigvee_{i \in I} [\underline{S}_{i}] = \bigwedge\{[\underline{T}]| [\underline{S}_{i}] \leq [\underline{T}]$ $\forall i \in I\} $. D{\"o}ring and Cannon also showed the following, which will be crucial for the purposes of this note, \\

	\textbf{Theorem 2.6:} $E$ and $L$ are isomorphic as complete lattices. In particular, the maps $g: E \rightarrow L$ and $f: L \rightarrow E$ defined by $g([\underline{S}]) = \varepsilon (\underline{S})$ and $f(a) = [\underline{\delta(a)}]$ are an inverse pair of complete lattice isomorphisms. \\

	\section{$\mathbb{R}^{(Sub_{cl}(\underline{\Omega}))} \sim SA(H)?$}

	The first thing to note is that since $Sub_{cl}(\underline{\Omega})$ is a complete Heyting algebra, $V^{Sub_{cl}(\underline{\Omega})}$ will be a `Heyting valued model' of intuitionistic set theory. Now, to begin the desired proof, we need to show that given any element $u \in V^{Sub_{cl}(\underline{\Omega})}$ of our model such that $\|$$u$ is a Dedekind real number$\| = \top$\footnote{Note that we need to specify that $u$ is a Dedekind real, since there are models of intuitionistic set theory where Dedekind and Cauchy reals are distinct. See \cite{HVM} for details.}, there corresponds a unique spectral family of projections on the the original Hilbert space, and vice-versa. Towards this end, let $u$ be such an element. Then we can consider the truth values $P_{q} = \| \hat{q} \in u \|$ ($q \in \mathbb{Q}$) with properties\footnote{These properties all follow immediately from the definition of Dedekind cuts.}\\

	\hspace{2cm}  (i) $\bigwedge _{r \in \mathbb{Q}}$ $P_{r}$ = $\bot$ \\

\hspace{2cm}  (ii) $\bigvee _{r \in \mathbb{Q}}$ $P_{r}$ = $\top$ \\

\hspace{2cm}  (iii) For any rational number r, $\bigwedge _{s > r}$ $P_{s} = P_{r}$ \\

	We can go on to define $E_{\lambda} = \bigwedge _{q > \lambda} P_{q}$ for $\lambda \in \mathbb{R}$, and establish the following\\

\hspace{2cm}  (i') $\bigwedge _{\lambda}$ $E_{\lambda}$ = $\bot $\\

\hspace{2cm}  (ii') $\bigvee _{\lambda}$ $E_{\lambda}$ = $\top$ \\

\hspace{2cm}  (iii') $\bigwedge _{\mu > \lambda}$ $E_{\mu} = E_{\lambda}$ \\

	Now, let $\phi(u) = \{E_{\lambda}|  \lambda \in \mathbb{R}\}$. We can then compose $\phi$ with the map $[.]$ that takes each $E_{\lambda}$ to its equivalence class. Thus, we have $([.] \circ \phi) (u) = \{ [E_{\lambda}]| \lambda \in \mathbb{R}\}$. Of the following properties,  (i'') and (iii'') are immediate from (i')-(iii') and the way that the lattice operations are defined on $E$. However, (ii'') is not guaranteed to hold because of the fact that $\varepsilon$ does not preserve arbitrary joins. \\

	\hspace{2cm}  (i'') $\bigwedge _{\lambda}$ $[E_{\lambda}]$ = $[\bot]$ \\

\hspace{2cm}  (ii'') $\bigvee _{\lambda}$ $[E_{\lambda}]$ = $[\top]$ \\

\hspace{2cm}  (iii'') $\bigwedge _{\mu > \lambda}$ $[E_{\mu}] = [E_{\lambda}]$ \\

	By theorem 2.6, we can use the isomorphism $g: E \rightarrow P(H)$ to obtain the set $\{ g([E_{\lambda}]) | \lambda \in \mathbb{R}\} \subseteq P(H)$. Since $g$ is an isomorphism, we immediately obtain the following properties, where we write $G_{\lambda}$ for $g([E_{\lambda}])$ \\

	\hspace{2cm}  (i''') $\bigwedge _{\lambda}$ $G_{\lambda}$ = $0$ \\

\hspace{2cm}  (iii''') $\bigwedge _{\mu > \lambda}$ $G_{\mu} = G_{\lambda}$ \\

	Unfortunately, $\{G_{\lambda} | \lambda \in \mathbb{R}\}$ is not quite a spectral family of projections, since it is not guaranteed to satisfy\\

\hspace{2cm}  (ii''') $\bigvee _{\lambda}$ $G_{\lambda}$ = $1$ \\

	The fact that $\{ g([E_{\lambda}]) | \lambda \in \mathbb{R}\}$ satisfies (i''') and (iii'''), but not (ii''') means that $\{ g([E_{\lambda}]) | \lambda \in \mathbb{R}\}$ is what is known as a `weakly right continuous spectral family'\footnote{This terminology comes from \cite{DD}, where these kinds of families are introduced and studied in the context of defining observables in Topos Quantum Theory. }. Intuitively, these families correspond to `self adjoint' operators that have positive infinity in their spectrum (let's call these operators `weakly self adjoint'). Thus, we have defined a map ($g \circ [.] \circ \phi$) that takes an element of $\mathbb{R}^{(Sub_{cl}(\underline{\Omega}))}$ and returns a weakly self-adjoint operator on the relevant Hilbert space. Let's call this entire composite map `$G$'. 

	Conversely, given a bounded\footnote{The requirement that $A$ be bounded is necessary for technical reasons that will be discussed shortly.} self-adjoint operator $A$ on the relevant Hilbert space, we want to find a corresponding element of $\mathbb{R}^{(Sub_{cl}(\underline{\Omega}))}$. We begin by making another application of the spectral theorem to obtain a bounded, left continuous spectral family $\{A_{\lambda}| \lambda \in \mathbb{R}\} \subseteq P(H)$ of projections. Next, we use $g$'s inverse isomorphism $f: P(H) \rightarrow E$ to obtain the set $\{f(A_{\lambda})| \lambda \in \mathbb{R}\} \subseteq E$. Obviously, we have\\

	\hspace{2cm}  (1) $\bigwedge _{\lambda}$ $F_{\lambda}$ = $[\bot]$ \\

\hspace{2cm}  (2) $\bigvee _{\lambda}$ $F_{\lambda}$ = $[\top]$ \\

\hspace{2cm}  (3) $\bigvee _{\mu < \lambda}$ $F_{\mu} = F_{\lambda}$ \\

	where $F_{\lambda} = f(A_{\lambda})$. \\

 At this stage, we need to define a new map from $E$ to $Sub_{cl}(\underline{\Sigma})$ with some special properties. This is achieved by the following (easily verified) lemma. \\

	\textbf{Lemma 3.1:} Define $h: E \rightarrow Sub_{cl}(\underline{\Omega})$ by $h([\underline{S}]) = \delta(\epsilon(\underline{S}))$ ($h$ is obviously well defined). Then $h$ satisfies  (i) $\varepsilon(h([\underline{S}])) = \varepsilon(\underline{S})$, $\forall [\underline{S}] \in E$ (ii) $h$ preserves joins (iii) $h$ is injective.\\

	So, we can now consider the set $\{h(F_{\lambda})| \lambda \in \mathbb{R}\} \subseteq  Sub_{cl}(\underline{\Sigma})$. We will write $H_{\lambda}$ for $h(F_{\lambda})$. Since $h$ preserves joins, it follows immediately that (2') and (3') below hold \\

	\hspace{2cm}  (1') $\bigwedge _{\lambda}$ $H_{\lambda}$ = $\bot$ \\

\hspace{2cm}  (2') $\bigvee _{\lambda}$ $H_{\lambda}$ = $\top$ \\

\hspace{2cm}  (3') $\bigvee _{\mu < \lambda}$ $H_{\mu} = H_{\lambda}$ \\

	To prove (1'), we need to use the fact that $A$ is bounded from below. Specifically, we know that there must be some $\lambda \in \mathbb{R}$ such that $A_{\lambda} = 0$. Since $f$ is an isomorphism, $F_{\lambda} = [\bot]$. It is easy to see that $H_{\lambda} = \bot$ has to hold, and so (1') is guaranteed. So $\{H_{\lambda}| \lambda \in \mathbb{R}\}$ is a set of truth values in $ Sub_{cl}(\underline{\Sigma})$ satisfying (1') - (3'). This allows us to define a corresponding real number in $V^{ Sub_{cl}(\underline{\Sigma})}$ in the same way as in Takeuti's original proof, i.e. we just define $u \in V^{ Sub_{cl}(\underline{\Sigma})}$ by $dom(u) = \{ \hat{q} \: : q \in \mathbb{Q} \}$ and $u(\hat{q}) = H_{q}$. It easily follows that $\| u$ is a Dedekind real number $\| = \top$. Thus, we have defined a map from $BSA(H) \subseteq SA(H)$ (the set of all bounded self adjoint operators) to $\mathbb{R}^{(Sub_{cl}(\underline{\Omega}))}$, as desired. We call the whole composite map `$H$'. The next step in our proof is to show that $G$ and $H$, as defined above, are both injective.

	We already know that $f: P(H) \rightarrow E$, being an isomorphism, is injective. Furthermore, we know that $h$ is injective, by lemma 3.1. So, given two distinct (bounded, left continuous) spectral families of projections, $\{A_{\lambda}| \lambda \in \mathbb{R}\}$, $\{B_{\lambda} | \lambda \in \mathbb{R}\}$, there will be $\lambda \in \mathbb{R}$ such that $A_{\lambda} \neq B_{\lambda}$. By injectivity of $f$, $F^{A}_{\lambda} \neq F^{B}_{\lambda}$. By injectivity of $h$, $H^{A}_{\lambda} \neq H^{B}_{\lambda}$. This guarantees that $\|u^{A} = u^{B}\| \neq \top$, i.e. $H$ is injective. So we have defined an injection of $BSA(H)$ into $\mathbb{R}^{(Sub_{cl}(\underline{\Omega}))}$.

	Conversely, if we have $u, v \in  \mathbb{R}^{(Sub_{cl}(\underline{\Omega}))}$ such that $\| u = v\| \neq \top$, then there will be $\lambda \in \mathbb{R}$ such that $E^{u}_{\lambda} \neq E^{v}_{\lambda}$. Now, we want to show that there will also exist some $\lambda \in \mathbb{R}$ such that $[E^{u}_{\lambda}] \neq [E^{v}_{\lambda}]$. In order to do this, we need to show that there exists $q \in \mathbb{Q}$ such that $[\|\hat{q} \in u\|] \neq [\|\hat{q} \in v\|]$. The most natural way to do this is to use the following template for a proof.\\

	\textbf{First Attempt At A Proof:} Let $u, v \in \mathbb{R}^{(Sub_{cl}(\underline{\Omega}))}$ with $\| u = v\| \neq \top$, and assume for contradiction that $\forall q \in \mathbb{Q}$, $[\|\hat{q} \in u\|] = [\|\hat{q} \in v\|]$. Then, since $\|u, v$ are Dedekind real numbers$\| = \top$, we have \\

	$\|u = v\| = \bigwedge_{q \in \mathbb{Q}} \|\hat{q} \in u \leftrightarrow \hat{q} \in v\|$\\

	$\hspace{1.3cm} = \bigwedge_{q \in \mathbb{Q}} \|\hat{q} \in u \rightarrow \hat{q} \in v\| \wedge \|\hat{q} \in v \rightarrow \hat{q} \in u\|$\\

	$\hspace{1.3cm} = \bigwedge_{q \in \mathbb{Q}} \|\hat{q} \in u \rightarrow \hat{q} \in v\| \wedge  \bigwedge_{q \in \mathbb{Q}} \|\hat{q} \in v \rightarrow \hat{q} \in u\|$\\
	
	Now, let $q \in \mathbb{Q}$. Then,\\

	 $[\|\hat{q} \in u \rightarrow \hat{q} \in v\|] = [\|\hat{q} \in u\| \Rightarrow \|\hat{q} \in v\|] $\\

	\hspace{2.7cm} $= [\neg \|\hat{q} \in u\| \vee \| \hat{q} \in v\|]$\\

	\hspace{2.7cm} $ \geq [\neg\|\hat{q} \in u\|] \vee [\| \hat{q} \in v\|]$\\

	\hspace{2.7cm} $ = [\neg\|\hat{q} \in u\|] \vee [\| \hat{q} \in u\|]$\\

	\hspace{2.7cm} $ = [\top]$\\

	At this stage, we need another lemma. \\

\textbf{Lemma 3.2:} $[\top] = \{\top\}$, i.e. $\varepsilon (\underline{S}) = \varepsilon(\top)$ implies $\underline{S} = \top$. \\

	\textbf{Proof:} Suppose that $\varepsilon (\underline{S}) = \varepsilon(\top)$. So $\varepsilon (\underline{S}) = 1$. So $\top = \underline{\delta(1)} = \underline{\delta(\varepsilon(S))} \leq \underline{S}$\\

	So, for any $q \in \mathbb{Q}$, $[\|\hat{q} \in u \rightarrow \hat{q} \in v\|] = [\top]$. So, by lemma 3.2, for any $q \in \mathbb{Q}$, $\|\hat{q} \in u \rightarrow \hat{q} \in v\| = \top$. So $\bigwedge_{q \in \mathbb{Q}} \|\hat{q} \in u \rightarrow \hat{q} \in v\| = \top$. The same argument shows that $\bigwedge_{q \in \mathbb{Q}} \|\hat{q} \in v \rightarrow \hat{q} \in u\| = \top$, and hence that $\| u = v\| = \top$, contradicting our assumption. So there must be some $q \in \mathbb{Q}$ such that $[\|\hat{q} \in u\|] \neq [\|\hat{q} \in v\|]$, as desired. 

	The above proof works perfectly except for one unjustified step. Specifically, the final step in the chain of equalities leading from $[\|\hat{q} \in u \rightarrow \hat{q} \in v\|]$ to $[\top]$ is illegitimate. It is not generally the case that for any clopen subobject $\underline{S}$, $[\neg \underline{S}] \vee [\underline{S}] = [\top]$. Indeed, since the clopen subobjects of the spectral presheaf only form a Heyting algebra, not a Boolean algebra, it is not even generally true that $\neg \underline{S} \vee \underline{S} = \top$, and this is weaker than the desired condition. 

	Now, the most natural solution to this problem is to use the fact that the clopen subobjects actually form a complete \emph{bi-Heyting algebra}, i.e. they also form a complete co-Heyting algebra\footnote{For technical details, see \cite{D1}} that comes equipped with a paraconsistent negation $\sim$ satisfying $\sim\underline{S} \vee \underline{S} = \top$ (but violating $\sim \underline{S} \wedge \underline{S} = \bot)$. This suggests that it might be possible to run the above proof, using the paraconsistent rather than the intuitionistic logical structure of the clopen subobjects. We simply go through the proof and replace every occurrence of the Heyting negation and implication with the corresponding co-Heyting operations. However, this still does not solve our problem since it turns out that it is not generally true that for any clopen subobject, $[\underline{S}] \vee [\sim \underline{S}] = [\top]$. Thus, of the two logical structures of the clopen subobjects that have been studied in the literature, neither appears to be able to facilitate the proof that is necessary to establish the desired connection between quantum set theory and topos quantum logic. In order to rectify this situation, we need to introduce a third, entirely new logical structure that has not previously been considered.

	\section{* - The Third Quantum Negation }

	\textbf{Def 4.1:} Given $\underline{S} \in Sub_{cl}(\underline{\Omega})$, define $\underline{S}^{*} = \underline{\delta(\varepsilon(S)^{\bot})}$, i.e. $\underline{S}^{*}$ is the daseinisation of the orthocomplement of $\varepsilon(\underline{S})$ ($\bot$ denotes the orthocomplement of $P(H)$. Recall that $\bot$ satisfies (a) $a \vee a^{\bot} = 1$, (b) $a \wedge a^{\bot} = 0$, (c) $a \leq b$ implies $a^{\bot} \geq b^{\bot}$, (d) $a^{\bot \bot} = a$, (e) $(a \wedge b)^{\bot} = a^{\bot} \vee b^{\bot}$, (f) $(a \vee b)^{\bot} = a^{\bot} \wedge b^{\bot}$). \\

	\textbf{Theorem 4.2:} The $*$ operation has the following properties, \\

	(i) $\underline{S} \vee \underline{S^{*}} = \top$\\

	(ii) $\underline{S}^{**} \leq \underline{S}$\\

	(iii) $\underline{S}^{***} = \underline{S}^{*}$\\

	(iv) $\underline{S} \wedge \underline{S}^{*} \geq \bot$\\

	(v) $(\underline{S} \wedge \underline{T})^{*} = \underline{S}^{*} \vee \underline{T}^{*}$\\

	(vi) $(\underline{S} \vee \underline{T})^{*} \leq \underline{S}^{*} \wedge \underline{T}^{*}$\\

	(vii) $\varepsilon(\underline{S}) \vee \varepsilon(\underline{S}^{*}) = 1$\\

	(viii) $\varepsilon(\underline{S}) \wedge \varepsilon(\underline{S}^{*}) = 0$\\

	(ix) $\underline{S} \leq \underline{T}$ implies $\underline{S}^{*} \geq \underline{T}^{*}$\\

		These properties show that $*$ plays the role of a particular kind (the details of which will be considered in the next section) of paraconsistent negation. For current purposes, the key property is (vii), which guarantees that for any clopen subobject $\underline{S}$, $[\underline{S}^{*}] \vee [\underline{S}] = [\top]$. Thus, if we let $\underline{S} \Rightarrow \underline{T}$ denote $\underline{S}^{*} \vee \underline{T}$ and replace each occurrence of the Heyting negation with $*$, then our attempted proof in the previous section works perfectly. 

	One nice feature of $*$ is that it allows us to extend the isomorphism between $E$ and $P(H)$ to include the orthocomplement operation on $P(H)$. Specifically, if we define $[\underline{S}]^{*} = [\underline{S}^{*}]$ (this is obviously well defined, since $\varepsilon(\underline{S}) = \varepsilon(\underline{T})$ clearly implies $\underline{S}^{*} = \underline{T}^{*})$, then we have $g([\underline{S}]^{*}) = \varepsilon(\underline{S}^{*}) = \varepsilon(\underline{\delta(\varepsilon(S)^{\bot})}) = \varepsilon(\underline{S})^{\bot} = g([\underline{S}])^{\bot}$. Also, $f(a^{\bot}) = [\underline{\delta(a^{\bot})}]$ and $f(a)^{*} = [\underline{\delta(a)}]^{*} = [\underline{\delta(a)}^{*}]$. Now, $\varepsilon(f(a^{\bot})) = a^{\bot}$ and $\varepsilon(\underline{\delta(a)}^{*}) = \varepsilon(\underline{\delta(a)})^{\bot} = a^{\bot}$. This proves that $f(a)^{*} = f(a^{\bot})$. So, under the negation operation $*$, $E$ is isomorphic to $L$ equipped with its orthocomplement. So, we have now shown that, by using *, we can guarantee that for any elements $u, v \in \mathbb{R}^{(Sub_{cl}(\underline{\Omega}))}$ with $\| u = v\| \neq \top$, there will exist some $\lambda \in \mathbb{R}$ such that $[E^{u}_{\lambda}] \neq [E^{v}_{\lambda}]$. Then, since g is an isomorphism, we will have $g([E^{u}_{\lambda}]) \neq g([E^{v}_{\lambda}])$, i.e. $G^{u}_{\lambda} \neq G^{v}_{\lambda}$, i.e. $G$ is injective, as desired. 

	\section{Paraconsistent Set Theory}

	In the preceding section, it was shown that if we consider $Sub_{cl}(\underline{\Omega})$ not as a Heyting or co-Heyting algebra, but rather as equipped with the new paraconsistent negation $*$, then our attempted proof of the injectivity of $G$ goes through smoothly. However, the situation here is more complicated than it seems. For, since we are not treating $Sub_{cl}(\underline{\Omega})$ as a Heyting algebra,  $V^{Sub_{cl}(\underline{\Omega})}$ is no longer guaranteed to be a model of intuitionistic set theory, and it is certainly not going to model full $ZFC$. Indeed, the only set theory that our structure will be able to model would be a set theory built over a suitable paraconsistent logic. 

	The study of paraconsistent set theories is still a relatively new field that has not yet been greatly explored. Furthermore, there has not yet been any attempt in the literature to build the kind of algebraic valued models described in this chapter for any of the theories that have been developed in this context. 

	One recent attempt to build a set theory over a paraconsistent logic will be important for our purposes. This is Zach Weber's paraconsistent set theory, which we will refer to as `$PST$'. We will now give a brief introduction to the relevant technical details of $PST$ (for a full technical development, see \cite{PST}). 

	First of all, we need to describe the underlying first order logic of $PST$. This is characterised by the following axioms and rules \\

	\textbf{Axioms:} All instances of the following schemata are theorems, \\

	 (1) $\Phi \rightarrow \Phi$\\

	(2a) $\Phi \wedge \Psi \rightarrow \Phi$\\

	(2b) $\Phi \wedge \Psi \rightarrow \Psi$\\

	(3) $\Phi \wedge (\Psi \vee \Theta) \rightarrow (\Phi \wedge \Psi) \vee (\Phi \wedge \Theta) $\\

	(4) $(\Phi \rightarrow \Psi) \wedge (\Psi \rightarrow \Theta) \rightarrow (\Phi \rightarrow \Theta)$\\

	(5) $(\Phi \rightarrow \Psi) \wedge (\Phi \rightarrow \Theta) \rightarrow (\Phi \rightarrow \Psi \wedge \Theta)$\\

	(6) $(\Phi \rightarrow \neg \Psi) \rightarrow (\Psi \rightarrow \neg \Phi)$\\

	(7) $\neg \neg \Psi \rightarrow \Psi$\\

	(8) $\Phi \vee \neg \Phi$\\

	(9) $\forall x \Phi \rightarrow \Phi (y/x)$\\

	(10) $ \forall x (\Phi \rightarrow \Psi) \rightarrow (\Phi \rightarrow \forall x \Psi)$\\

	(11) $ \forall x (\Phi \vee \Psi) \rightarrow (\Phi \vee \forall x \Psi)$\\

	\textbf{Rules:} We also assume the following rules of inference, \\

	 (1) $\Phi, \Psi \vdash \Phi \wedge \Psi$\\

	(2) $\Phi, \Phi \rightarrow \Psi \vdash \Psi$\\

	(3) $\Phi, \neg \Psi \vdash \neg (\Phi \rightarrow \Psi)$\\

	(4) $\Phi \rightarrow \Psi, \Theta \rightarrow \Delta \vdash (\Psi \rightarrow \Theta) \rightarrow (\Phi \rightarrow \Delta$)\\

	(5) $\Phi \rightarrow \forall x \Phi$\\

	(6) $x = y \vdash \Phi(x) \rightarrow \Phi(y)$\\

	\textbf{Theorem 5.1:} $Sub_{cl}(\underline{\Omega})$, equipped with $*$ and the corresponding implication connective ($\underline{S} \Rightarrow \underline{T} = \underline{S}^{*} \vee \underline{T}$) is an algebraic model of  all of the axioms of the logic of $PST$, and all of the rules except for rule (3)\footnote{This logic is actually well known in the paraconsistent logic community. It is known as `dialectical logic with quantifiers' (see \cite{Rout}).}. \\

	We omit the proof here in the interest of concision, but the theorem follows mainly from theorem 4.2. The importance of theorem 5.1 is that $V^{Sub_{cl}(\underline{\Omega})}$ will be a structure whose internal logic is very closely related to the underlying logic of $PST$. This suggests that once the model theory of $PST$ has been properly developed, $V^{Sub_{cl}(\underline{\Omega})}$ it should be possible to treat $V^{Sub_{cl}(\underline{\Omega})}$ as one of its models.

	Now that we've become acquainted with the logic underlying PST, it is time to get introduced to the theory itself. One of the major advantages of PST, that stems from using a paraconsistent logic, is that the axiomatisation of the theory is far simpler than that of other well known set theories. Specifically, the theory consists of the following two axioms, \\

	\textbf{Axioms of PST:} \\

	(1) Abstraction: $x \in \{z| \phi(z)\} \leftrightarrow \phi(x)$.\\

	(2) Extensionality: $\forall z(z \in x \leftrightarrow z \in y) \leftrightarrow x = y$\\

	The abstraction schema is completely unrestricted, i.e. it tells us that for any formula $\phi$, we can form the set of all sets satisfying $\phi$, even if the set thus defined appears free in $\phi$. Of course, as usual, this leads to the familiar set-theoretic paradoxes like Burali-Forte and the Russell paradox. However, in $PST$, these contradictions can be tolerated without the system becoming trivial (for a proof of the non-triviality of $PST$, see \cite{Brady}) because the underlying logic is paraconsistent. Indeed, these paradoxes become basic theorems in $PST$. Specifically, the Burrali-Forte paradox becomes a theorem that tells us that the set of all ordinals (which we can form using the unrestricted abstraction scheme) both is and is not a member of itself. 

	Now, in a recent series of papers, Weber and others have shown how all the standard set-theoretic machinery (e.g basic algebra of sets, cardinal and ordinal arithmetic etc) from $ZFC$ can be constructed in $PST$. It has also been shown that $PST$ provides a new perspective on many deep and challenging problems in $ZFC$. For example, it has been shown \cite{PST} that in $PST$, the continuum hypothesis is provably false, while the existence of measurable cardinals is a theorem. Work has also been done on developing real analysis in the context of paraconsistent logics like that underlying $PST$. Again, the basic construction of Cauchy and Dedekind reals and the proof of their most important properties has all been recovered in this setting, and we know that $PST$ has the existence of the Dedekind reals as a theorem.

	\section{Conclusions and Further Work}

	In summary, we have seen that the structure $V^{Sub_{cl}(\underline{\Omega})}$, when equipped with the new paraconsistent negation $*$, looks like a very interesting potential model for $PST$. Specifically, assuming that we can legitimately treat $V^{Sub_{cl}(\underline{\Omega})}$ as a model of $PST$, we have established the existence of injections $G:\mathbb{R}^{(Sub_{cl}(\underline{\Omega}))} \rightarrow WSA(H)$ (where $WSA(H)$ is the set of all weakly self adjoint operators on $H$) and $H: BSA(H) \rightarrow \mathbb{R}^{(Sub_{cl}(\underline{\Omega}))}$. This leaves us with three main directions for further work, \\

	\begin{itemize}

	\item{Develop theorem transfer schemes for $V^{Sub_{cl}(\underline{\Omega})}$ similar to those obtained for $V^{(P(H))}$ in \cite{Oza1}}

	\item{Further develop our understanding of the cardinality of $\mathbb{R}^{(Sub_{cl}(\underline{\Omega}))}$, with an ultimate aim of recreating a Takeuti style bijection with $SA(H)$}

	\item{Study the relationship between $V^{Sub_{cl}(\underline{\Omega})}$ and the presheaf toposes of TQT. In particular, study how much of the internal machinery of TQT can be recreated `inside' $V^{Sub_{cl}(\underline{\Omega})}$}

	\item{Study how $V^{Sub_{cl}(\underline{\Omega})}$ can be used to translate concepts and results between TQT and QST, thereby expanding the expressive power and physical significance of these formalisms.}

	\end{itemize}

	\section{Acknowledgements}

	Many thanks to the QPL referees for their helpful comments. Special thanks also to Joerg Brendle, Andreas D{\"o}ring and Richard Pettigrew for all their patient help and guidance. The author was supported by an AHRC doctoral scholarship during the writing of this paper. 

\bibliographystyle{eptcs}
\bibliography{generic}

\end{document}